\documentclass[letterpaper, 10 pt, conference]{./ieeeconf/ieeeconf}   

\IEEEoverridecommandlockouts                              

\usepackage{graphics} 
\usepackage{graphicx} 
\usepackage{epsfig} 

\usepackage{mathptmx} 

\usepackage{amsmath} 
\usepackage{amssymb}  
\usepackage[caption=false]{subfig}
\usepackage[draft]{hyperref}

\makeatletter

\newcommand{\Rmnum}[1]{\expandafter\@slowromancap\romannumeral #1@}
\makeatother

\title{\LARGE \bf
	Change of Optimal Values: A Pre-calculated Metric
}

\author{Fang Bai
\thanks{Fang Bai is with Center for Autonomous Systems, University of Technology Sydney (UTS-CAS), Sydney, Australia. Email:
	fang.bai@yahoo.com, fang.bai@student.uts.edu.au}%
}

\newtheorem{remark}{\textbf{Remark}}

\begin{document}

\maketitle
\thispagestyle{empty}
\pagestyle{empty}

\begin{abstract}

A variety of optimization problems takes the form of a minimum norm optimization.
In this paper, we study the change of optimal values between two incrementally constructed least norm optimization problems, with new measurements included in the second one.
We prove an exact equation to calculate the change of optimal values in the linear least norm optimization problem.
With the result in this paper,
the change of the optimal values can be pre-calculated as a metric to guide online decision makings,
without solving the second optimization problem as long the solution and covariance of the first optimization problem are available.
The result can be extended to linear least distance optimization problems, and nonlinear least distance optimization with (nonlinear) equality constraints through linearizations.
This derivation in this paper provides a theoretically sound explanation to the empirical observations shown in \cite{bai2018robust}.
As an additional contribution, we propose another optimization problem, i.e. aligning two trajectories at given poses, to further demonstrate how to use the metric.
The accuracy of the metric is validated with numerical examples,
which is quite satisfactory in general (see the experiments in \cite{bai2018robust} as well),
unless in some extremely adverse scenarios.
Last but not least, calculating the optimal value by the proposed metric is at least one magnitude faster than solving the corresponding optimization problems directly.

\end{abstract}


\section{Introduction}

Many robotic applications rely on metrics to make decisions.
A metric is usually formed as a pre-calculated scalar which captures certain aspects of the underlying problem, for example the entropy in the information theory,
or the cost in an optimization problem.

The most widely used metrics in robotics, at least in estimation tasks, are attributed to the information theoretic metrics \cite{mackay2003information}.
This is a broad class, and we mention several of them for completeness:
(1) Fisher Information Matrix (FIM).
If the processing noise is Gaussian, the metric always boils down to the evaluation of a proper covariance matrix with respect to its eigen value, trace or determinant \cite{carrillo2012comparison}\cite{khosoussi2018reliable}.
(2) Mutual Information.
This metric is quite popular which has been extensively used in robotics,
for example, the active control/planning strategy \cite{whaite1997autonomous}\cite{feder1999adaptive}\cite{bourgault2002information}\cite{davison2005active}\cite{vidal2010action},
graph pruning \cite{kretzschmar2012information}\cite{carlevaris2014generic},
online sparse pose-graph construction \cite{ila2010information}\cite{johannsson2013temporally},
path planning on a pose graph \cite{valencia2013planning}
and so forth.
(3) Kullback-Leibler Divergence (KL Divergence).
While the first two metrics deal merely with uncertainty, KL divergence can incorporate the mean and covariance of two Gaussian distributions together.
This is quite similar to the optimal value in the optimization context that we will discuss later, but they are obviously in different tracks.
See \cite{wang2013kullback} for an example of graph pruning using KL divergence.

Besides the information theoretic metrics, in the context of optimization,
the optimal value is essentially a metric as well.
However, the optimal value was considered dilatory in the past,
because the optimal value can only be obtained after solving the underlying problem.
This is meaningless because in most applications, we often encounter a bundle of candidate problems and a metric is required to tell us which problem is worth solving.
Obviously, if solving the problem is compulsory, the optimal value is unsuitable for tasks of this kind.
As a result, the optimal value is basically used for evaluating the optimality of the solution \cite{olson2009evaluating}\cite{hu2010evaluation}.
Some authors managed to use the optimal value \cite{latif2013robust}\cite{graham2015robust} to detect outliers by using an ``optimize-revoke" strategy.
However, this strategy is quite expensive as many undesirable subproblems are formed and solved in the process.

Recently, for an incremental optimization problem, Bai et al. \cite{bai2018robust} empirically show that the optimal value can be approximated prior to solving the problem.
In \cite{bai2018robust}, we provided an equation for the optimal value, but were unable to prove it exactly.
Later in \cite{bai2018predicting}, the authors studied incremental least squares optimization, and proved that the optimal value can be pre-calculated exactly for the linear case.
The theory extends to nonlinear cases through linearization.
However, there is still no clear explanation for the result on the equality constrained optimization in \cite{bai2018robust}.

In this paper, we study the least norm optimization, which is the counterpart of least squares.
Analogous to \cite{bai2018predicting}, we first derive the result for the linear case which is exact, and then extend it to more general problems like least distance optimization and nonlinear scenarios.
The analysis in this paper eventually provides a theoretical sound explanation for the empirical observation in \cite{bai2018robust}.
Moreover, the theoretical result in this paper is actually quite general, which can be applied to other similar problems easily.

Overall, this paper includes the following contributions:
\begin{enumerate}
	\item The change of optimal values in minimum norm optimization (Section \ref{section: Results on Least Norm Optimization}), and its extension to linear least distance optimization (Section \ref{section: Extension to Least Distance Optimization}).
	\item An extension to equality constrained nonlinear least distance optimization on manifold (Section \ref{section: Extension to Nonlinear Cases with an Example on Manifold}), which explains why the change of optimal values can be pre-calculated as empirically shown in the paper \cite{bai2018robust}.
	\item An example on the cost of aligning trajectories to illustrate the usage of the metric (Section \ref{section: Cost of Aligning Two Trajectories}), with numerical experiments in Section \ref{section: Numerical Experiments}.
\end{enumerate}

\section{Results on Least Norm Optimization}
\label{section: Results on Least Norm Optimization}

\subsection{Problem Statement and Standard Form}

Let us consider a minimum norm optimization problem
\begin{equation}
	\label{minimal norm optimzation: full stacked measurement}
	\min\  \mathbf{x}^T \mathbf{x} 
	\qquad
	\mathbf{s.t.}\quad
	\mathbf{A} \mathbf{x} = \mathbf{b}
\end{equation}
where the constraints are classified into two parts:
\begin{equation*}
\mathbf{A} = 
\begin{bmatrix}
\mathbf{A}_1 \\ \mathbf{A}_2
\end{bmatrix}
,
\qquad
\mathbf{b} = 
\begin{bmatrix}
\mathbf{b}_1 \\ \mathbf{b}_2
\end{bmatrix}
.
\end{equation*}
This is particularly the case for temporally incremental optimization problems.
The constraint
$
\mathbf{A}_1 \mathbf{x} = \mathbf{b}_1
$
represents those that are already available at a specific time point,
while the constraint
$
\mathbf{A}_2 \mathbf{x} = \mathbf{b}_2
$
can only be obtained after that time point.
Let us assume both $\mathbf{A}_1$ and $\mathbf{A}_2$ have full row rank.

The optimization problem established using the previously available measurement explicitly writes
\begin{equation}
\label{minimal norm optimzation: partial measurement, previous}
\min\  \mathbf{x}^T \mathbf{x} 
\qquad
\mathbf{s.t.}\quad
\mathbf{A}_1 \mathbf{x} = \mathbf{b}_1
.
\end{equation}
Now at some point, we have solved the problem (\ref{minimal norm optimzation: partial measurement, previous}), but not yet problem (\ref{minimal norm optimzation: full stacked measurement}).
The question is that: Can we tell something useful about problem (\ref{minimal norm optimzation: full stacked measurement}) from the results of problem (\ref{minimal norm optimzation: partial measurement, previous})?

In this section, we establish the equation on how to calculate the optimal value (of the objective function) of problem (\ref{minimal norm optimzation: full stacked measurement}), without solving it directly.

\subsection{Classical Results on Minimum Norm Optimization}

To proceed, let us recall some classical results of minimum norm optimization, taking problem
(\ref{minimal norm optimzation: full stacked measurement}) as an example.

The optimal solution $\mathbf{x}_{opt}$ is
\begin{equation*}
\mathbf{x}_{opt} 
=
\mathbf{A}^{\dagger} \mathbf{b}
=
\mathbf{A}^T (\mathbf{A} \mathbf{A}^T)^{-1} \mathbf{b}
.
\end{equation*}

The covariance (with Gaussian assumption) of the optimal solution, $\mathbb{C}\mathrm{ov} (\mathbf{x}_{opt} )$, is the projection matrix to the Null space of $\mathbf{A}$, which writes
\begin{equation*}
	\mathbb{C}\mathrm{ov} (\mathbf{x}_{opt} )
	=
	\mathbf{P}_{(N(\mathbf{A}))}
	= \mathbb{I} - \mathbf{A}^T (\mathbf{A} \mathbf{A}^T)^{-1} \mathbf{A}
	.
\end{equation*}

The cost at the optimal solution, termed the optimal value (of the objective function), i.e. $\mathbf{x}_{opt}^T \mathbf{x}_{opt}$, writes
\begin{equation*}
	f_{opt} =
	\mathbf{b}^T 
	(\mathbf{A} \mathbf{A}^T)^{-1}
	\mathbf{b}
	.
\end{equation*}

We will use the Schur complement \cite{meyer2000matrix} to decompose the below block matrix 
\begin{equation}
\label{equation: schur complement to decompose matrix}
	\begin{bmatrix}
	\mathbf{U} & \mathbf{P} \\ \mathbf{Q} & \mathbf{V}
	\end{bmatrix}
	=
	\begin{bmatrix}
	\mathbb{I} & \mathbb{O} \\ \mathbf{Q}\mathbf{U}^{-1} & \mathbb{I}
	\end{bmatrix}
	\begin{bmatrix}
	\mathbf{U} & \mathbb{O} \\ \mathbb{O} & 
	\mathbf{V} - \mathbf{Q} \mathbf{U}^{-1} \mathbf{P}
	\end{bmatrix}
	\begin{bmatrix}
	\mathbb{I} & \mathbf{U}^{-1}\mathbf{P}  \\\mathbb{O}   & \mathbb{I}
	\end{bmatrix}
	.
\end{equation}
Eq. (\ref{equation: schur complement to decompose matrix}) stands as long as $\mathbf{A}$ is invertible.

\subsection{Derivation of the Main Result}

The derivation consists of a manipulation of matrices via linear algebra laws.
Basically, we will expand the optimal value of problem (\ref{minimal norm optimzation: full stacked measurement}), i.e.,
$
\mathbf{b}^T 
(\mathbf{A} \mathbf{A}^T)^{-1}
\mathbf{b}
$,
and examine its components.

The two matrices below are useful to express the Schur complement of $\mathbf{A} \mathbf{A}^T$.
Let us define two invertible matrices:
\begin{equation}
\mathbf{W} = 
\mathbf{A}_2 
\mathbb{C}\mathrm{ov} (\mathbf{x}_1^{\star} )    \mathbf{A}_2^T
= 
\mathbf{A}_2 [\mathbb{I} - \mathbf{A}_1^T (\mathbf{A}_1 \mathbf{A}_1^T)^{-1} \mathbf{A}_1] \mathbf{A}_2^T
\end{equation}
\begin{equation}
\mathbf{H} = 
	\begin{bmatrix}
	\mathbb{I} &  \mathbb{O}  \\
	\mathbf{A}_2 \mathbf{A}_1^T (\mathbf{A}_1 \mathbf{A}_1^T)^{-1} &
	\mathbb{I}
	\end{bmatrix}
	=
	\begin{bmatrix}
	\mathbb{I} &  \mathbb{O}  \\
	\mathbf{A}_2 \mathbf{A}_1^{\dagger} &
	\mathbb{I}
	\end{bmatrix}
	.
\end{equation}

Let us apply Schur complement (\ref{equation: schur complement to decompose matrix}) to the matrix $\mathbf{A} \mathbf{A}^T$. The result writes
\begin{equation*}
\begin{aligned}
	\mathbf{A} 	\mathbf{A}^T
	= &
	\begin{bmatrix}
	\mathbf{A}_1 \mathbf{A}_1^T & \mathbf{A}_1 \mathbf{A}_2^T \\
	\mathbf{A}_2 \mathbf{A}_1^T &
	\mathbf{A}_2 \mathbf{A}_2^T
	\end{bmatrix}
	\\[3pt] = &
	\mathbf{H}
	\begin{bmatrix}
	\mathbf{A}_1 \mathbf{A}_1^T  &  \mathbb{O}  \\
	\mathbb{O} &
	\mathbf{A}_2 \mathbf{A}_2^T - \mathbf{A}_2 \mathbf{A}_1^T (\mathbf{A}_1 \mathbf{A}_1^T)^{-1} \mathbf{A}_1 \mathbf{A}_2^T
	\end{bmatrix}
	\mathbf{H}^T
	\\[3pt] = &
	\mathbf{H}
	\begin{bmatrix}
	\mathbf{A}_1 \mathbf{A}_1^T  &  \mathbb{O}  \\
	\mathbb{O} &
	\mathbf{W}
	\end{bmatrix}
	\mathbf{H}^T
	.
\end{aligned}
\end{equation*}

Noting that $\mathbf{H}^{-1} \mathbf{b}$ can be written as
\begin{equation*}
\mathbf{H}^{-1} \mathbf{b}
=
\begin{bmatrix}
\mathbb{I} &  \mathbb{O}  \\
- \mathbf{A}_2 \mathbf{A}_1^{\dagger} &
\mathbb{I}
\end{bmatrix}
\begin{bmatrix}
\mathbf{b}_1 \\ \mathbf{b}_2
\end{bmatrix}
=
\begin{bmatrix}
\mathbf{b}_1 \\ 
- \mathbf{A}_2 \mathbf{A}_1^{\dagger} \mathbf{b}_1 + \mathbf{b}_2
\end{bmatrix},
\end{equation*}
then we can finally expand
$
\mathbf{b}^T 
(\mathbf{A} \mathbf{A}^T)^{-1}
\mathbf{b}
$
by a straightforward matrix calculation as follow
\begin{equation*}
\begin{aligned}
& f^{\star\star} = 
	 \mathbf{b}^T 
	(\mathbf{A} \mathbf{A}^T)^{-1}
	\mathbf{b}
	\\[3pt] = &
	\mathbf{b}^{T} 
	\mathbf{H}^{-T}
	\begin{bmatrix}
	\mathbf{A}_1 \mathbf{A}_1^T  &  \mathbb{O}  \\
	\mathbb{O} &
	\mathbf{W}
	\end{bmatrix}^{-1}
	\mathbf{H}^{-1}
	\mathbf{b}	
	\\[3pt] = &
	(\mathbf{H}^{-1} \mathbf{b})^T
	\begin{bmatrix}
	(\mathbf{A}_1 \mathbf{A}_1^T)^{-1}  &  \mathbb{O}  \\
	\mathbb{O} &
	\mathbf{W}^{-1}
	\end{bmatrix}
	\mathbf{H}^{-1} \mathbf{b}	
	\\[3pt] = &
	\begin{bmatrix}
	\mathbf{b}_1 \\ 
	- \mathbf{A}_2 \mathbf{A}_1^{\dagger} \mathbf{b}_1 + \mathbf{b}_2
	\end{bmatrix}^T
	\begin{bmatrix}
	(\mathbf{A}_1 \mathbf{A}_1^T)^{-1} &  \mathbb{O}  \\
	\mathbb{O} &
	\mathbf{W}^{-1}
	\end{bmatrix}
	\begin{bmatrix}
	\mathbf{b}_1 \\ 
	- \mathbf{A}_2 \mathbf{A}_1^{\dagger} \mathbf{b}_1 + \mathbf{b}_2
	\end{bmatrix}	
	\\[4pt] = &
	\mathbf{b}_1^T 	(\mathbf{A}_1 \mathbf{A}_1^T)^{-1} \mathbf{b}_1
	+
	(\mathbf{A}_2 \mathbf{A}_1^{\dagger} \mathbf{b}_1 - \mathbf{b}_2)^T
	\mathbf{W}^{-1}
	(\mathbf{A}_2 \mathbf{A}_1^{\dagger} \mathbf{b}_1 - \mathbf{b}_2)
	\\[4pt] = &
	f^{\star} + 
	\underbrace{(\mathbf{A}_2 \mathbf{x}_1^{\star} - \mathbf{b}_2)^T
		[\mathbf{A}_2 
		\mathbb{C}\mathrm{ov} (\mathbf{x}_1^{\star} )    \mathbf{A}_2^T]^{-1}
		(\mathbf{A}_2 \mathbf{x}_1^{\star} - \mathbf{b}_2)}_{\Delta f}
	.
\end{aligned}
\end{equation*}

\textbf{Conclusion:}
Therefore the \textbf{change of optimal values},
$
\Delta f = f^{\star\star} - f^{\star}
$,
from problem (\ref{minimal norm optimzation: full stacked measurement}) to problem (\ref{minimal norm optimzation: partial measurement, previous})
can be compactly written as
\begin{equation}
\label{equation: predicting objective function change in minimal norm optimization}
\Delta f 
=
(\mathbf{A}_2 \mathbf{x}_1^{\star} - \mathbf{b}_2)^T
[\mathbf{A}_2 
\mathbb{C}\mathrm{ov} (\mathbf{x}_1^{\star} )    \mathbf{A}_2^T]^{-1}
(\mathbf{A}_2 \mathbf{x}_1^{\star} - \mathbf{b}_2)
.
\end{equation}
In (\ref{equation: predicting objective function change in minimal norm optimization}), the unknown piece of information,
$\mathbf{x}_1^{\star}$ and $\mathbb{C}\mathrm{ov} (\mathbf{x}_1^{\star} )$,
can be obtained by solving the problem (\ref{minimal norm optimzation: partial measurement, previous}).
In other words, given $\mathbf{x}_1^{\star}$ and its covariance $\mathbb{C}\mathrm{ov} (\mathbf{x}_1^{\star} )$, we can calculate $\Delta f $ through Eq. (\ref{equation: predicting objective function change in minimal norm optimization}), then we calculate the the optimal value of problem (\ref{minimal norm optimzation: full stacked measurement}) simply by $f^{\star\star} = f^{\star} + \Delta f$, without solving it!

\section{Extension to Least Distance Optimization}
\label{section: Extension to Least Distance Optimization}

In this section,
we consider the problem of a least distance optimization, whose standard form writes
\begin{equation}
\label{Eq: least distance optimization, full measurement matrix}
\min\  (\mathbf{H}\mathbf{x} - \mathbf{h})^T \boldsymbol{\Sigma}^{-1} (\mathbf{H}\mathbf{x} - \mathbf{h})
\qquad
\mathbf{s.t.}\quad
\mathbf{A} \mathbf{x} = \mathbf{b}
.
\end{equation}
In (\ref{Eq: least distance optimization, full measurement matrix}),
we assume $\mathbf{H}$, $\boldsymbol{\Sigma}$ to be invertible.
In this case,
by letting $\mathbf{y} = \boldsymbol{\Sigma}^{-\frac{1}{2}}(\mathbf{H}\mathbf{x} - \mathbf{h})$, i.e.
$
\mathbf{x} = \mathbf{H}^{-1}\boldsymbol{\Sigma}^{\frac{1}{2}}\mathbf{y} + \mathbf{H}^{-1}\mathbf{h}
$,
(\ref{Eq: least distance optimization, full measurement matrix}) can be converted to a least norm optimization
\begin{equation}
\min\  \mathbf{y}^T \mathbf{y}
\qquad
\mathbf{s.t.}\quad
\mathbf{A} \mathbf{H}^{-1} \boldsymbol{\Sigma}^{\frac{1}{2}} \mathbf{y} = \mathbf{b} - \mathbf{A} \mathbf{H}^{-1} \mathbf{h}
.
\end{equation}
Therefore we can extend the conclusion on the least norm optimization to
(\ref{Eq: least distance optimization, full measurement matrix}).

As before, we assume $\mathbf{A}$ and $\mathbf{b}$ are obtained in two different phases comprising
\begin{equation*}
\mathbf{A} = 
\begin{bmatrix}
\mathbf{A}_1 \\ \mathbf{A}_2
\end{bmatrix}
,
\qquad
\mathbf{b} = 
\begin{bmatrix}
\mathbf{b}_1 \\ \mathbf{b}_2
\end{bmatrix}
\end{equation*}
with $\mathbf{A}_1$, $\mathbf{A}_2$ of full row rank.
Let us write the optimization problem using the measurement from the first phase only as
\begin{equation}
\label{Eq: least distance optimization, partial measurement matrix}
\min\  (\mathbf{H}\mathbf{x} - \mathbf{h})^T \boldsymbol{\Sigma}^{-1} (\mathbf{H}\mathbf{x} - \mathbf{h})
\qquad
\mathbf{s.t.}\quad
\mathbf{A}_1 \mathbf{x} = \mathbf{b}_1
.
\end{equation}

Formally, we would like to calculate the change of optimal values from the problem
(\ref{Eq: least distance optimization, partial measurement matrix}) to the problem
(\ref{Eq: least distance optimization, full measurement matrix}).

Let us write (\ref{Eq: least distance optimization, partial measurement matrix})
in the form of a least norm optimization
\begin{equation}
\min\  \mathbf{y}^T \mathbf{y}
\qquad
\mathbf{s.t.}\quad
\mathbf{A}_1 \mathbf{H}^{-1} \boldsymbol{\Sigma}^{\frac{1}{2}} \mathbf{y} = \mathbf{b}_1 - \mathbf{A} \mathbf{H}^{-1} \mathbf{h}
\end{equation}
and denote its optimal solution as $\mathbf{y}_1^{\star}$.
The covariance of $\mathbf{y}_1^{\star}$ is related with that of
$\mathbf{x}_1^{\star} = \mathbf{H}^{-1}\boldsymbol{\Sigma}^{\frac{1}{2}}\mathbf{y}_1^{\star} + \mathbf{H}^{-1}\mathbf{h}
$
by the linear transformation of Gaussian distribution
\begin{equation}
	\label{Constrained Least Squares: intermediate equation: minimal norm optimization, covariance}
	\mathbf{A}_2 \mathbf{H}^{-1} \boldsymbol{\Sigma}^{\frac{1}{2}}
	\mathbb{C}\mathrm{ov} (\mathbf{y}_{1}^{\star} ) 
	\boldsymbol{\Sigma}^{\frac{1}{2}} \mathbf{H}^{-T} \mathbf{A}_2^T
	= 
	\mathbf{A}_2
	\mathbb{C}\mathrm{ov} (\mathbf{x}_{1}^{\star} )
	\mathbf{A}_2^T
	,
\end{equation}
and considering
$\mathbf{y}_1^{\star} = \boldsymbol{\Sigma}^{-\frac{1}{2}}(\mathbf{H}\mathbf{x}_1^{\star} - \mathbf{h}_1)
$
we have
\begin{equation}
	\label{Constrained Least Squares: intermediate equation: minimal norm optimization, residual}
	\mathbf{A}_2 \mathbf{H}^{-1} \boldsymbol{\Sigma}^{\frac{1}{2}} \mathbf{y}_1^{\star} - (\mathbf{b}_2 - \mathbf{A}_2 \mathbf{H}^{-1} \mathbf{h}_2)
	= \mathbf{A}_2 \mathbf{x}_1^{\star} - \mathbf{b}_2
	.
\end{equation}

Inserting (\ref{Constrained Least Squares: intermediate equation: minimal norm optimization, covariance}) and (\ref{Constrained Least Squares: intermediate equation: minimal norm optimization, residual})
into the result of minimum norm optimization,
we immediately conclude that the change of optimal values writes
\begin{equation}
\label{equation: predicting objective function change in least distance optimization}
	\Delta f =
	(\mathbf{A}_2 \mathbf{x}_1^{\star} - \mathbf{b}_2)^T
	[\mathbf{A}_2
	\mathbb{C}\mathrm{ov} (\mathbf{x}_{1}^{\star} )
	\mathbf{A}_2^T]^{-1}
	(\mathbf{A}_2 \mathbf{x}_1^{\star} - \mathbf{b}_2)
	.
\end{equation}

At last, we can write down a formula of
$
\mathbb{C}\mathrm{ov}(\mathbf{x}_{1}^{\star} )
$
by considering together
\begin{equation*}
\mathbb{C}\mathrm{ov} (\mathbf{y}_{1}^{\star} )
= \mathbb{I} - \boldsymbol{\Sigma}^{\frac{1}{2}} \mathbf{H}^{-T}\mathbf{A}_1^T (\mathbf{A}_1 \mathbf{H}^{-1} \boldsymbol{\Sigma} \mathbf{H}^{-T} \mathbf{A}_1^T)^{-1} \mathbf{A}_1 \mathbf{H}^{-1} \boldsymbol{\Sigma}^{\frac{1}{2}}
\end{equation*}
\begin{equation*}
\mathbb{C}\mathrm{ov} (\mathbf{x}_{1}^{\star} )
= 
\mathbf{H}^{-1}\boldsymbol{\Sigma}^{\frac{1}{2}}
\mathbb{C}\mathrm{ov} (\mathbf{y}_{1}^{\star} )
\boldsymbol{\Sigma}^{\frac{1}{2}}\mathbf{H}^{-T}
\end{equation*}
to reach its explicit form
\begin{equation}
\label{eq: covariance of least distance optimization}
\mathbb{C}\mathrm{ov} (\mathbf{x}_{1}^{\star} )
= 
\mathbf{Q} - \mathbf{Q} \mathbf{A}_1^T (\mathbf{A}_1 \mathbf{Q} \mathbf{A}_1^T)^{-1} \mathbf{A}_1 \mathbf{Q} 
\end{equation}
where $\mathbf{Q} = \mathbf{H}^{-1} \boldsymbol{\Sigma} \mathbf{H}^{-T}$.

To conclude, if we have known
$\mathbf{x}_{1}^{\star}$ and
$\mathbb{C}\mathrm{ov} (\mathbf{x}_{1}^{\star} )$
by solving the problem (\ref{Eq: least distance optimization, partial measurement matrix}),
then we can predict the change of optimal values between (\ref{Eq: least distance optimization, partial measurement matrix}) and (\ref{Eq: least distance optimization, full measurement matrix}) by
(\ref{equation: predicting objective function change in least distance optimization}) easily.
In the linear case, the result is exact.
For nonlinear cases, a linear approximation may apply around its current working point.

\section{Extension to Nonlinear Cases with an Example on Manifold}
\label{section: Extension to Nonlinear Cases with an Example on Manifold}

A nonlinear least distance optimization with equality constraints can be generalized as
\begin{equation*}
	\min\  \mathbb{D}\mathrm{ist} (\mathbf{x}, \tilde{\mathbf{x}})
	\qquad
	\mathbf{s.t.}\quad
	\mathbf{C}(\mathbf{x}) = \mathbf{0}
\end{equation*}
with $\tilde{\mathbf{x}}$ being the measurement data having the same dimension as $\mathbf{x}$,
and
$
\mathbb{D}\mathrm{ist} (\mathbf{x}, \tilde{\mathbf{x}})
$
being a scalar metric.

Usually, there are many ways to define a distance function
$
\mathbb{D}\mathrm{ist} (\mathbf{x}, \tilde{\mathbf{x}})
$.
Here we specifically focus on the Mahalanobis distance, which extends well to the Frobenius norm and $\ell_2$ norm based distance.
Explicitly, the problem we will discuss takes the form of
\begin{equation}
\label{Nonlinear Least Distance Optimization, Mahalanobis distance}
\min\  \lVert \mathbb{D}\mathrm{iff} (\mathbf{x}, \tilde{\mathbf{x}}) \rVert_{\Sigma}^2
\qquad
\mathbf{s.t.}\quad
\mathbf{C}(\mathbf{x}) = \mathbf{0}
\end{equation}
where $\mathbb{D}\mathrm{iff} (\mathbf{x}, \tilde{\mathbf{x}})$ is a vectorized difference between $\mathbf{x}$ and $\tilde{\mathbf{x}}$.
$\Sigma$ models our confidence on $\mathbb{D}\mathrm{iff} (\mathbf{x}, \tilde{\mathbf{x}})$ which is usually interpreted as the covariance of some Gaussian distribution.

A standard extension is to linearize (\ref{Nonlinear Least Distance Optimization, Mahalanobis distance})
to the form of
(\ref{Eq: least distance optimization, full measurement matrix}),
and then apply the result in linear cases,
as the paradigm found in Extended Kalman Filter (EKF) techniques.
In essence, the objective function in the linearized problem is an \textit{approximation} to that of its original nonlinear scenario, which means the result is not exact any more.
We mention EKF as an example because (\ref{equation: predicting objective function change in least distance optimization}) works incrementally, which shows some similarities to the Kalman Filter (KF).
\textit{For an incremental optimization problem, the linearization point will not change dramatically.}
Therefore it is reasonable to assume that the linearized Jacobian provides a valid approximation to the original nonlinear problem.

In what follows, we illustrate the idea of using linearization with an example on manifold.

Let $\mathbf{x}$ and $\tilde{\mathbf{x}}$ be two elements on the same manifold.
We use $\boxminus$ to describe their difference in the tangent space,
i.e. $\mathbb{D}\mathrm{iff} (\mathbf{x}, \tilde{\mathbf{x}}) = \mathbf{x} \boxminus \tilde{\mathbf{x}}$.
At a certain time point, we have an optimization problem, say \textbf{Phase\ \Rmnum{1}}:
\begin{equation*}
\begin{aligned}
  & \min\ \lVert \mathbf{x} \boxminus \tilde{\mathbf{x}} \rVert_{\Sigma}^2
\\
\mathbf{s.t.} \quad
& \mathbf{C}_1(\mathbf{x}) = \mathbf{0}
\\
\end{aligned}
\qquad(\mathbf{Phase\ \Rmnum{1}})
\end{equation*}
Let $\mathbf{x}_1^{\star}$ and $\mathbb{C}\mathrm{ov} (\mathbf{x}_1^{\star})$ be the solution and covariance of the optimization problem in \textbf{Phase \Rmnum{1}}.
Now if we obtain another constraint $\mathbf{C}_2(\mathbf{x}) = \mathbf{0}$,
which constitutes another optimization problem named \textbf{Phase \Rmnum{2}}:
\begin{equation*}
\begin{aligned}
 & \min\  \lVert \mathbf{x} \boxminus \tilde{\mathbf{x}} \rVert_{\Sigma}^2
\\
\mathbf{s.t.} \quad
& \mathbf{C}_1(\mathbf{x}) = \mathbf{0}
\\
& \mathbf{C}_2(\mathbf{x}) = \mathbf{0}
\\
\end{aligned}
\qquad(\mathbf{Phase\ \Rmnum{2}})
\end{equation*}
Then the change of optimal values from \textbf{Phase \Rmnum{1}} to \textbf{Phase \Rmnum{2}} can be approximated by
\begin{equation}
\label{equation: predicting objective function change in nonlinear least distance optimization}
\Delta f =
\mathbf{C}_2(\mathbf{x}_1^{\star})^T
[\mathbf{A}_2
\mathbb{C}\mathrm{ov} (\mathbf{x}_{1}^{\star} )
\mathbf{A}_2^T]^{-1}
\mathbf{C}_2(\mathbf{x}_1^{\star})
\end{equation}
with $\mathbf{A}_2$ given by
\begin{equation*}
\mathbf{A}_2 = - \frac{\partial \mathbf{C}_2 (\mathbf{x}_1^{\star} \boxplus \boldsymbol{\xi})}{\partial \boldsymbol{\xi}} \rvert_{\boldsymbol{\xi} = \mathbf{0}}
.
\end{equation*}
The notation $\boxplus$ is used to apply a perturbation $\xi$ in the tangent space of $\mathbf{x}_1^{\star}$, and return an element on manifold.
The covariance $\mathbb{C}\mathrm{ov} (\mathbf{x}_{1}^{\star} )$
is given in the tangent space of the solution $\mathbf{x}_{1}^{\star}$.
It is calculated by by Eq. (\ref{eq: covariance of least distance optimization}), by letting
\begin{equation*}
\mathbf{A}_1 = - \frac{\partial \mathbf{C}_1 (\mathbf{x}_1^{\star} \boxplus \boldsymbol{\xi})}{\partial \boldsymbol{\xi}} \rvert_{\boldsymbol{\xi} = \mathbf{0}}
\end{equation*}
\begin{equation*}
\mathbf{H} = - \frac{\partial \left( \left( \mathbf{x}_1^{\star} \boxplus \boldsymbol{\xi} \right) \boxminus \tilde{\mathbf{x}} \right)}{\partial \boldsymbol{\xi}} \rvert_{\boldsymbol{\xi} = \mathbf{0}}.
\end{equation*}

The derivation of Eq. (\ref{equation: predicting objective function change in nonlinear least distance optimization}) is quite straightforward via linearization. Let
\begin{equation*}
\mathbf{h} = \mathbf{x}_1^{\star} \boxminus \tilde{\mathbf{x}},
\qquad
\mathbf{b}_1 = \mathbf{C}_1( \mathbf{x}_1^{\star} ),
\qquad
\mathbf{b}_2 = \mathbf{C}_2( \mathbf{x}_1^{\star} )
.
\end{equation*}
At $\mathbf{x}_{1}^{\star}$, let us linearize the problem \textbf{Phase \Rmnum{1}} by a perturbation $\boldsymbol{\xi}$. The linearized problem writes
\begin{equation}
\label{eq: linearized least distance optimziation problem of Phase I}
\min\  (\mathbf{H}\boldsymbol{\xi} - \mathbf{h})^T \boldsymbol{\Sigma}^{-1} (\mathbf{H}\boldsymbol{\xi} - \mathbf{h})
\qquad
\mathbf{s.t.}\quad
\mathbf{A}_1 \boldsymbol{\xi} = \mathbf{b}_1
.
\end{equation}
Since $\mathbf{x}_{1}^{\star}$ is optimal for \textbf{Phase \Rmnum{1}}, the solution of the above linearized problem is $\xi_1^{\star} = 0$.
Now if we linearize the constraint $\mathbf{C}_2(\mathbf{x}) = \mathbf{0}$ to $\mathbf{A}_2 \boldsymbol{\xi} = \mathbf{b}_2$ at $\mathbf{x}_{1}^{\star}$
and add it in (\ref{eq: linearized least distance optimziation problem of Phase I}),
then by the result in the linear case, i.e. (\ref{equation: predicting objective function change in least distance optimization}),
\begin{equation*}
\begin{aligned}
	\Delta f = &
	 (\mathbf{A}_2 \boldsymbol{\xi}_1^{\star} - \mathbf{b}_2)^T
	[\mathbf{A}_2
	\mathbb{C}\mathrm{ov} (\boldsymbol{\xi}_{1}^{\star} )
	\mathbf{A}_2^T]^{-1}
	(\mathbf{A}_2 \boldsymbol{\xi}_1^{\star} - \mathbf{b}_2)
	\\[4pt] = &
	 \mathbf{b}_2^T
	[\mathbf{A}_2
	\mathbb{C}\mathrm{ov} (\boldsymbol{\xi}_{1}^{\star} )
	\mathbf{A}_2^T]^{-1}
	\mathbf{b}_2
	\\[4pt] = &
	\mathbf{C}_2( \mathbf{x}_1^{\star} )^T
	[\mathbf{A}_2
	\mathbb{C}\mathrm{ov} (\mathbf{x}_{1}^{\star} )
	\mathbf{A}_2^T]^{-1}
	\mathbf{C}_2( \mathbf{x}_1^{\star} )_2
\end{aligned}
\end{equation*}
which completes the proof.

\begin{remark}
Eq. (\ref{equation: predicting objective function change in nonlinear least distance optimization}) firstly appeared in the paper \cite{bai2018robust} (Eq. (15) in that paper). Unfortunately at the time,
the mechanism behind the equation was not completely clear.
So the analysis here also fills the gap in \cite{bai2018robust}.
\end{remark}

\section{Applications}

\subsection{Outlier Detection in Constrained Pose Graph Formulation}

The constrained pose graph formulation was firstly proposed in \cite{bai2018robust}, while the original idea on the feature based SLAM was published in \cite{bai2016incremental}.
In \cite{bai2018robust}, the change of optimal values was used as a metric \cite{bai2018robust} to detect outliers in an incremental scheme.
Although the mechanism was not clear in that paper, the experiments provided in \cite{bai2018robust} demonstrated the effectiveness of the change of optimal values as a metric to improve the robustness of the framework.
As an additional material, for the conventional least squares based pose graph formulation, the outlier detection based on the change of optimal values was published in \cite{bai2018predicting}.

\subsection{Cost of Aligning Two Trajectories}
\label{section: Cost of Aligning Two Trajectories}

\begin{figure*}[t]
	\centering
	{\includegraphics[width=.3\textwidth]{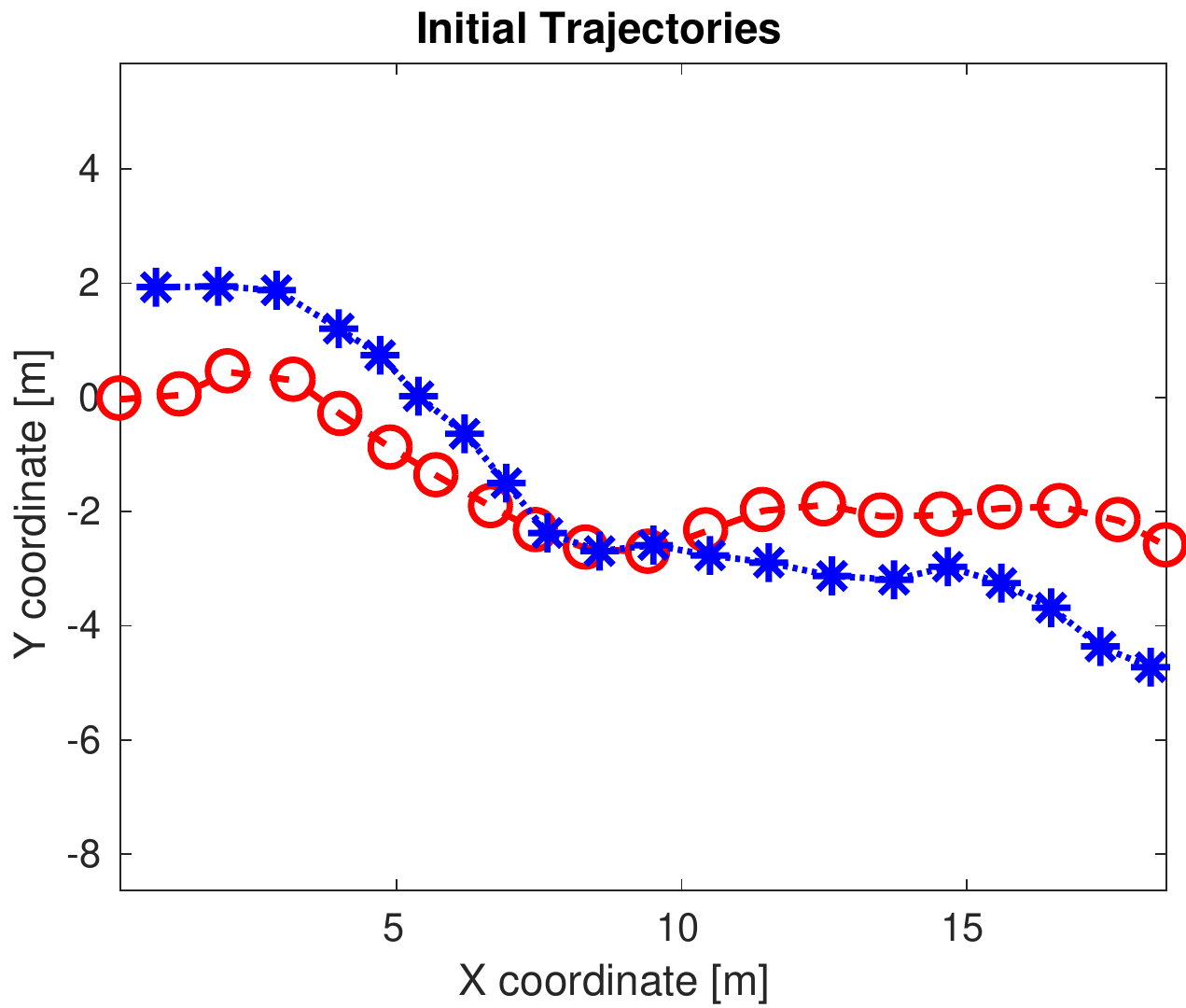}}
	{\includegraphics[width=.3\textwidth]{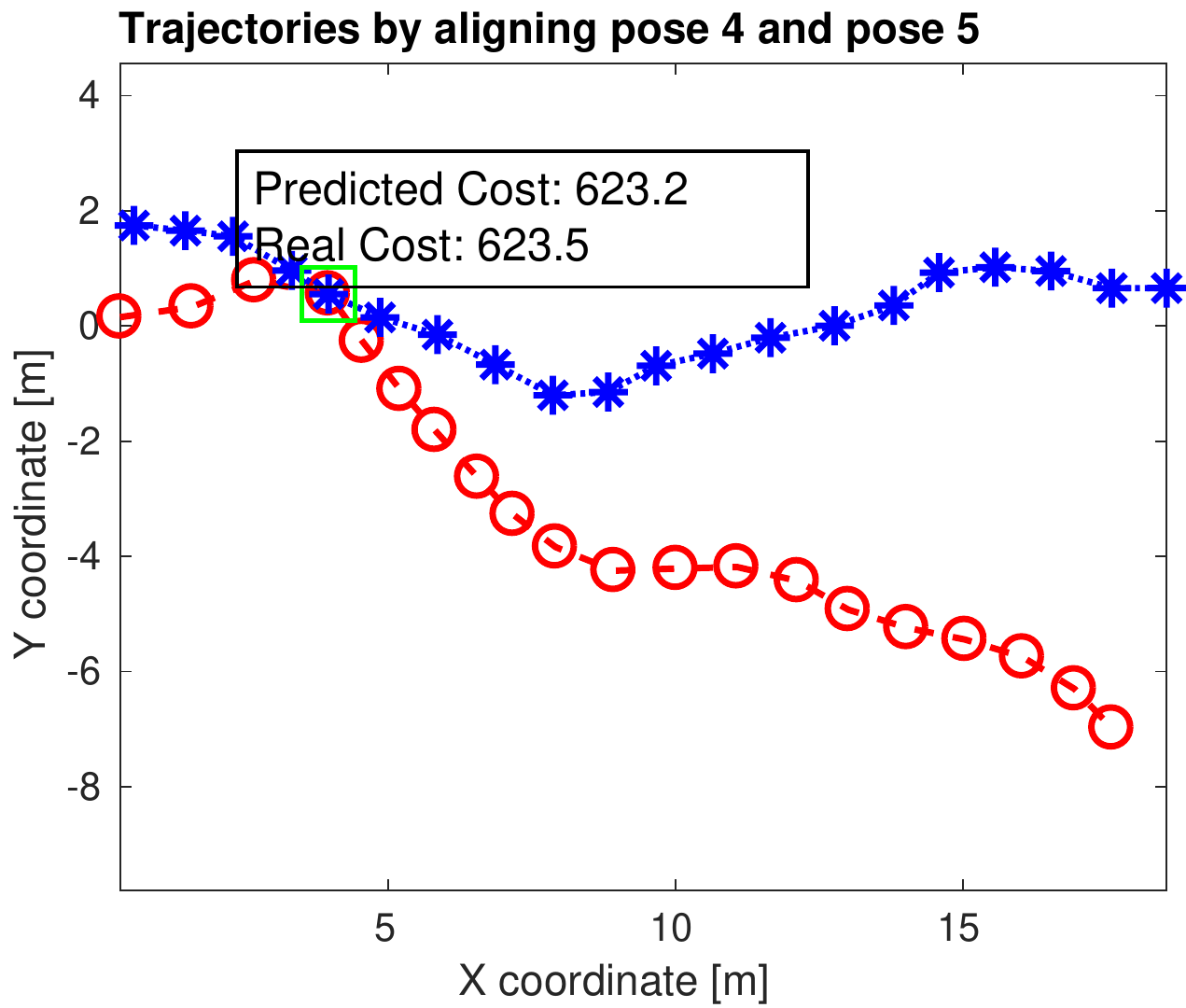}}
	{\includegraphics[width=.3\textwidth]{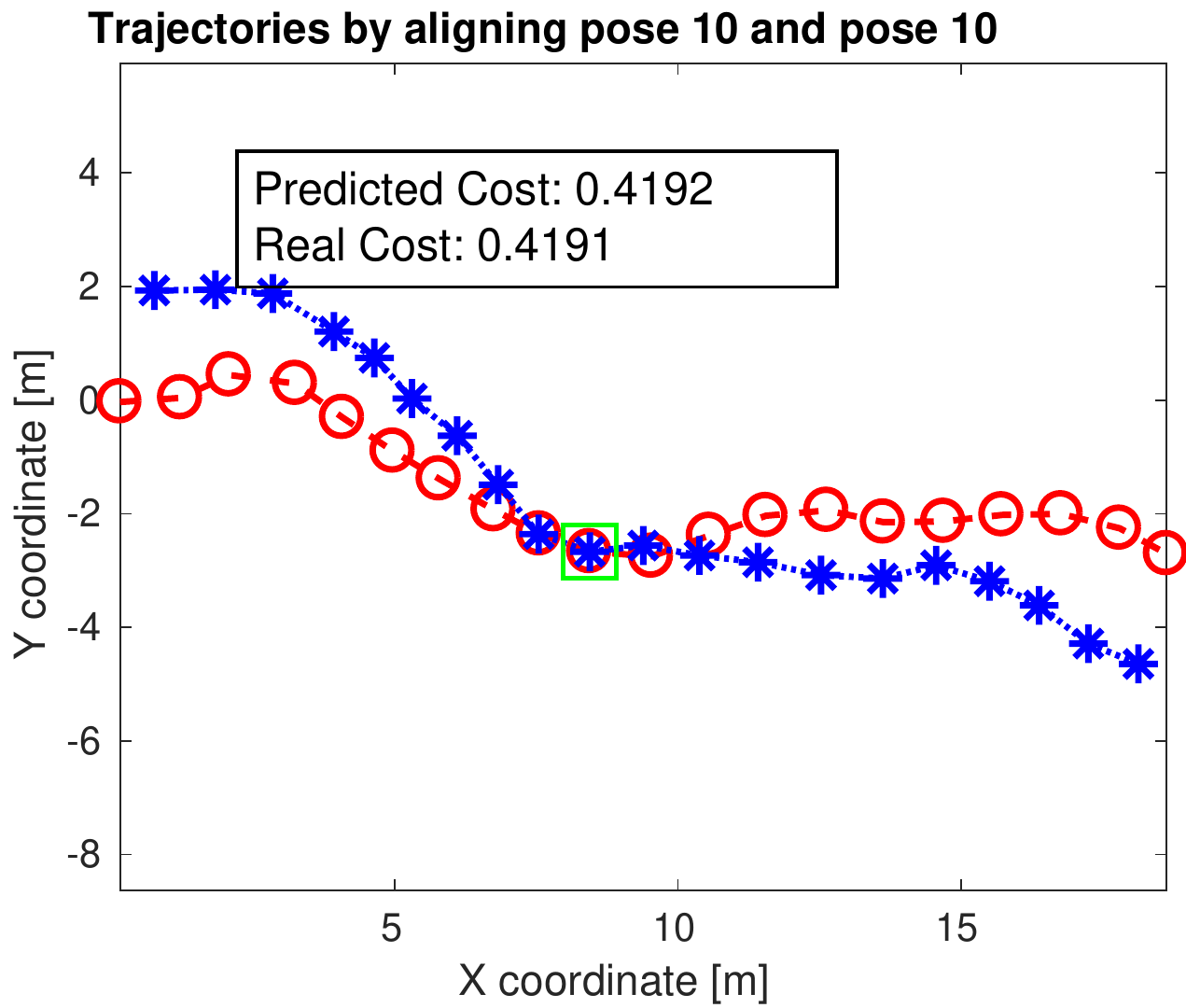}}
	{\includegraphics[width=.3\textwidth]{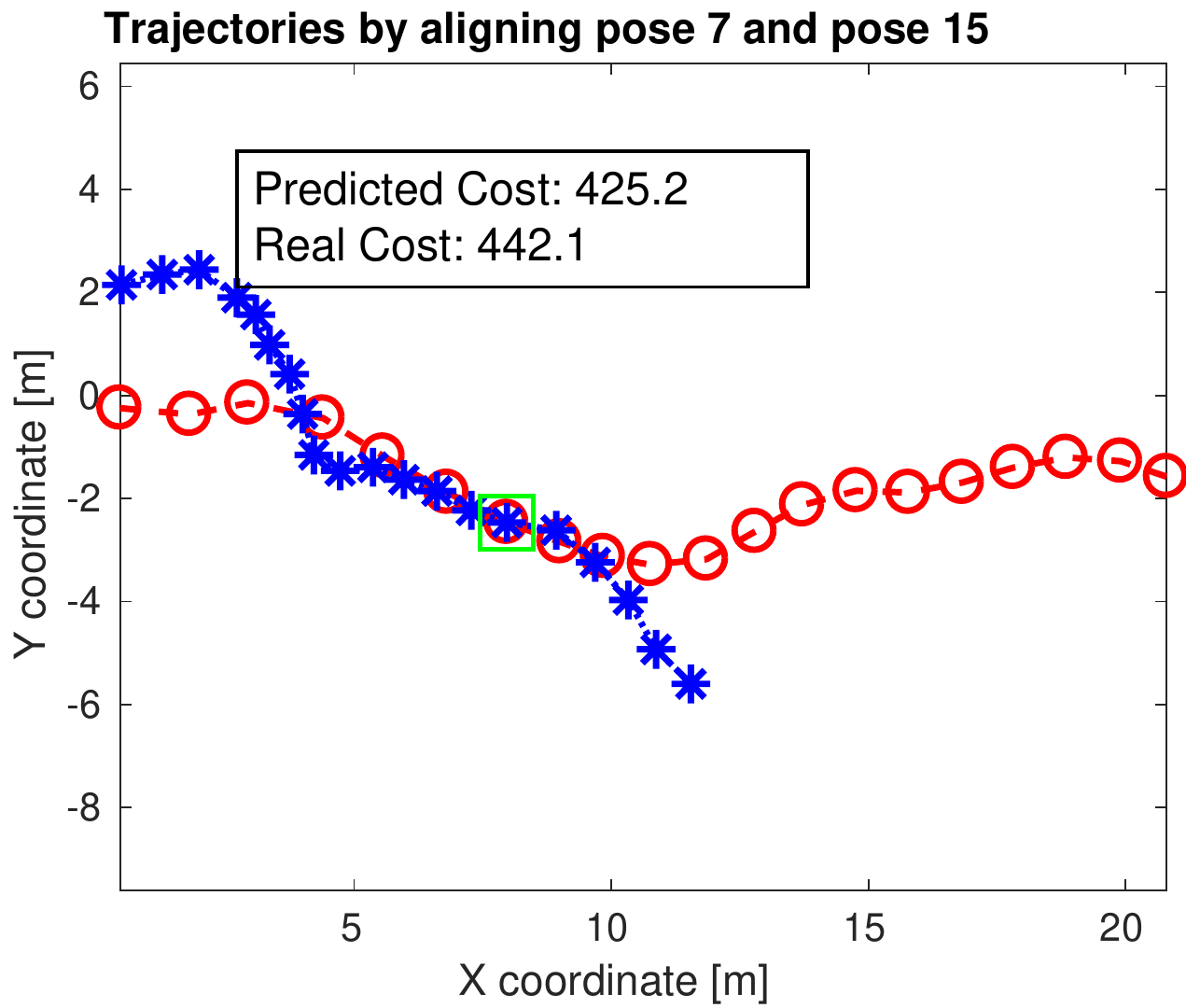}}
	{\includegraphics[width=.3\textwidth]{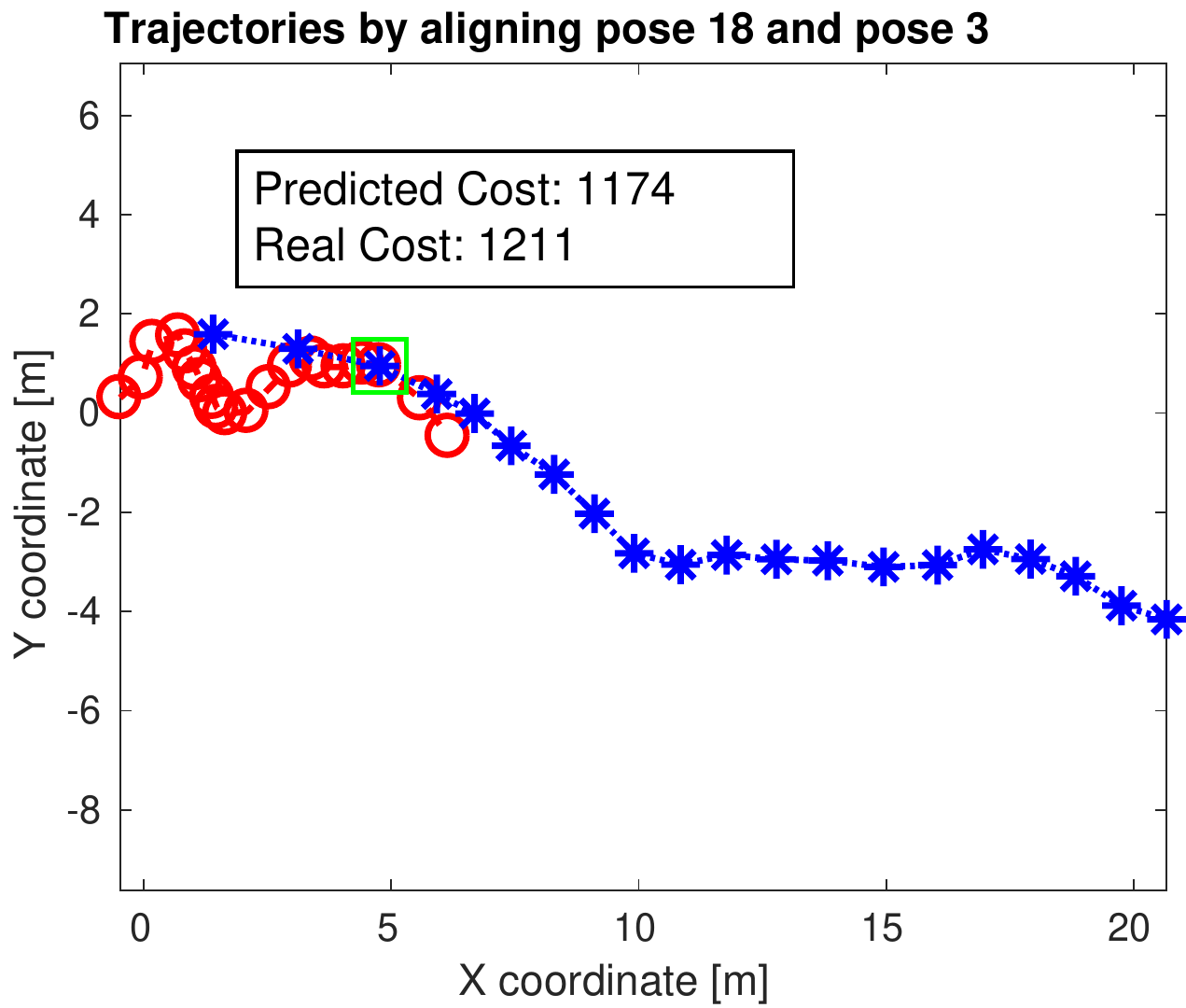}}
	{\includegraphics[width=.3\textwidth]{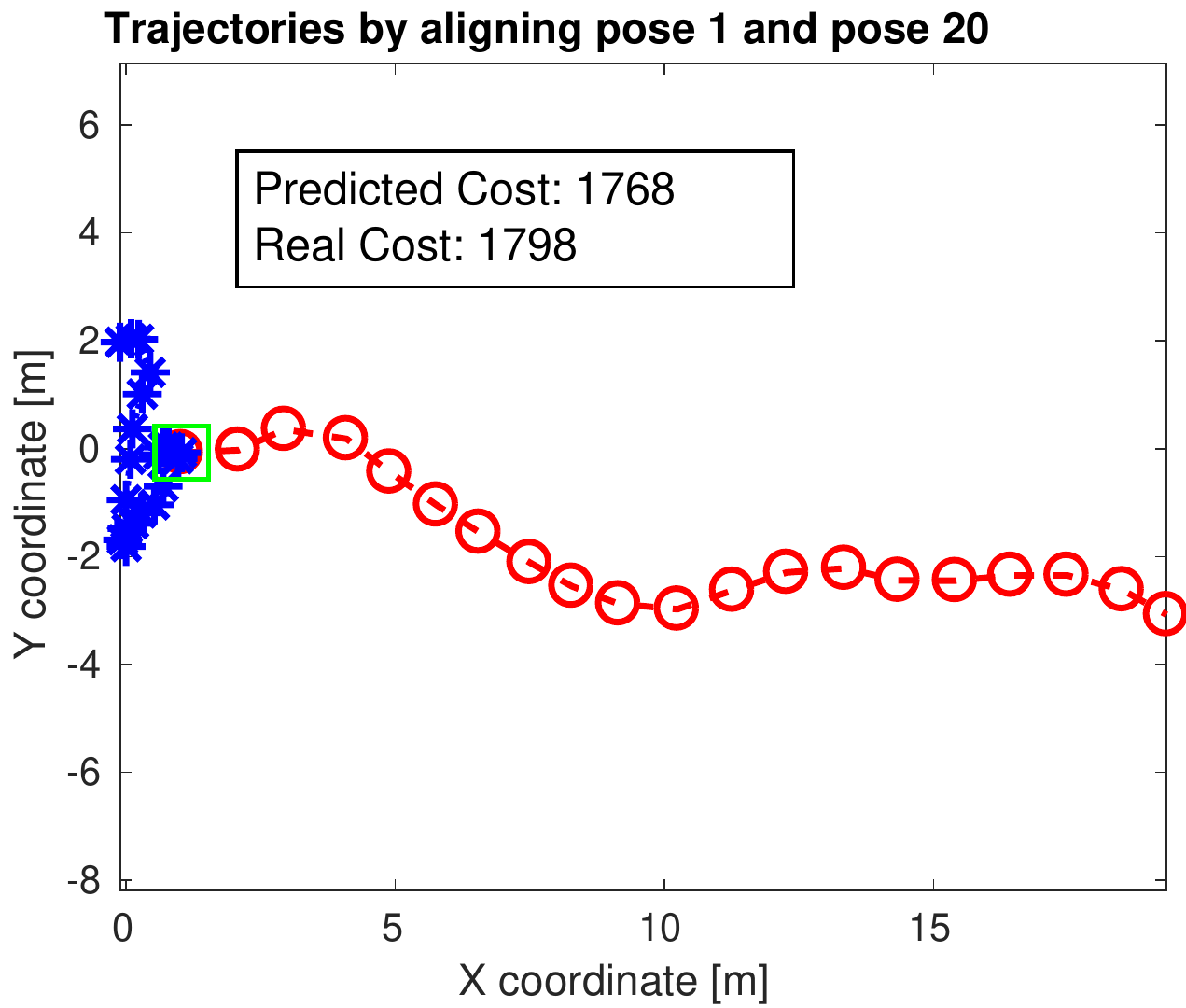}}		
	\caption{Examples of aligning two trajectories, plotted in red and blue respectively. For each aligned case, we show the resulting trajectory and the cost it takes. Both the predicted cost (by the OFC metric) and the real cost (by fully optimizing the problem) are reported for comparison.}
	\label{lb: Examples of aligning two trajectories, trajectory change and cost}
\end{figure*}

A trajectory is, in its explicit form, a collection of poses.
Let $\mathcal{T}$ be a trajectory, then it can be described by
$\mathcal{T} = \cup \mathbf{T}_i$,
$ \mathbf{T}_i \in \mathrm{SE}(d)$, $d = 2,3$.
However, in many applications, if one pose is moved, we want the impact to be passed to other poses as well.
Therefore, it is more desirable to exploit an alternative parameterization based on relative transformations (termed relative poses in robotics).
Without loss of generality, let us define a trajectory as a head pose, plus a chain of relative poses following the head pose.
Let $\mathcal{T}_A$ and $\mathcal{T}_B$ be two trajectories given by
\begin{equation*}
\mathcal{T}_A
\stackrel{\mathrm{def}}{=}
\{ {}^A\mathbf{T}_1 \} \cup
\{ {}^A\mathbf{T}_{i,i+1} \}_{i = 1, 2,\cdots N_A}
\end{equation*}
\begin{equation*}
\mathcal{T}_B
\stackrel{\mathrm{def}}{=} 
\{ {}^B\mathbf{T}_1 \} \cup
\{ {}^B\mathbf{T}_{j,j+1} \}_{j = 1, 2,\cdots N_B}
\end{equation*}
%
%
where ${}^A\mathbf{T}_1$ is the head pose of trajectory $\mathcal{T}_A$ and
$\{ {}^A\mathbf{T}_{i,i+1} \}_{i = 1, 2,\cdots N_A}$ is a set of relative poses following the head pose ${}^A\mathbf{T}_1$.
The same rule applies to trajectory $\mathcal{T}_B$.
In this parameterization, we obtain any arbitrary pose by chaining out the relative poses.
For example, in the trajectory $\mathcal{T}_A$, an arbitrary pose ${}^A\mathbf{T}_k$ can be written as
\begin{equation*}
	{}^A\mathbf{T}_k = {}^A\mathbf{T}_1  {}^A\mathbf{T}_{1,2}  \cdots  {}^A\mathbf{T}_{k-1,k}
	.
\end{equation*}

Now at some point, we know that the $l$-th pose ${}^A\mathbf{T}_l$ in the trajectory $\mathcal{T}_A$ and the $r$-th pose ${}^B\mathbf{T}_r$ in the trajectory $\mathcal{T}_B$ are the same pose, which results to a hard constraint
$
{}^A\mathbf{T}_l = {}^B\mathbf{T}_r
$
and an optimization problem in the form of
\begin{equation}
	\label{eq. optimization problem: cost of aliging trajectory gaps}
	\begin{aligned}
		\min\ &
	\underbrace{
			\lVert   {}^A\mathbf{T}_1   \boxminus  {}^A\tilde{\mathbf{T}}_1   \rVert_{{}^A\Sigma_1}^2
			 +
			\sum_{i=1}^{N_A} \lVert   {}^A\mathbf{T}_{i,i+1} \boxminus {}^A\tilde{\mathbf{T}}_{i,i+1}   \rVert_{{}^A\Sigma_{i,i+1}}^2	
		}_{\mathrm{Cost\ of\ trajectory\ } \mathcal{T}_A }
		\\  + &
	\underbrace{
			\lVert   {}^B\mathbf{T}_1  \boxminus  {}^B\tilde{\mathbf{T}}_1   \rVert_{{}^B\Sigma_1}^2		
			+
			\sum_{j=1}^{N_B} \lVert
			{}^B\mathbf{T}_{j,j+1}  \boxminus  {}^B\tilde{\mathbf{T}}_{j,j+1}  \rVert_{{}^B\Sigma_{j,j+1}}^2	
		}_{\mathrm{Cost\ of\ trajectory\ } \mathcal{T}_B }
		\\
		 \mathbf{s.t.} \quad
		 & {}^A\mathbf{T}_l = {}^B\mathbf{T}_r 		\\
		 &\ \mathrm{with}\ 
		  {}^A\mathbf{T}_l = {}^A\mathbf{T}_1 {}^A\mathbf{T}_{1,2} \cdots {}^A\mathbf{T}_{l-1,l}
		  \\ & \ \mathrm{\ \ \ \ \ }\ 	  
		  {}^B\mathbf{T}_r = {}^B\mathbf{T}_1 {}^B\mathbf{T}_{1,2} \cdots {}^B\mathbf{T}_{r-1,r}
	\end{aligned}
	.
\end{equation}
where we use the notations
${}^A\tilde{\mathbf{T}}_1$,
${}^B\tilde{\mathbf{T}}_1$,
${}^A\tilde{\mathbf{T}}_{i,i+1}$,
${}^B\tilde{\mathbf{T}}_{i,i+1}$
to describe the initial configuration of
$\mathcal{T}_A$ and $\mathcal{T}_B$
without the hard constraint
${}^A\mathbf{T}_l = {}^B\mathbf{T}_r$,
while
${}^A\mathbf{T}_1$,
${}^B\mathbf{T}_1$,
${}^A\mathbf{T}_{i,i+1}$,
${}^B\mathbf{T}_{i,i+1}$
are state variables to be estimated by imposing the hard constraint
${}^A\mathbf{T}_l = {}^B\mathbf{T}_r$.

\begin{remark}
	The problem described in (\ref{eq. optimization problem: cost of aliging trajectory gaps}) is essentially an extension to the problem ``trajectory bending" studied by Dubbelman et al. \cite{dubbelman2010efficient}\cite{dubbelman2015cop}.
	Please refer to \cite{dubbelman2010efficient}\cite{dubbelman2015cop} for a reference of its application backgrounds.
	Here in (\ref{eq. optimization problem: cost of aliging trajectory gaps}), we define the ``bending" as a more general constraint by aligning any arbitrary two poses together; while the original paper was to ``bend" a pose in a trajectory to a specific given value. However mathematically, I think this difference is marginal and one should extend to another easily.
	Nevertheless, it is worth noticing how seemingly different problems are essentially a description of kind of the same thing.
\end{remark}

\textbf{Problem:}
Due to process noise, there might exist many potential candidate pairs for ${}^A\mathbf{T}_l$ and ${}^B\mathbf{T}_r$.
The question is: What is the cost for each specific pair?
A brutal force answer is to solve the optimization problem in (\ref{eq. optimization problem: cost of aliging trajectory gaps}) for each pair.
However, alternatively we can use the theory described in this paper to predict the change of optimal values directly.

In (\ref{eq. optimization problem: cost of aliging trajectory gaps}),
let us choose $\boxminus$ as the logarithm mapping of $\mathrm{SE}(3)$ in the vector space.
To apply the result in Section \ref{section: Extension to Nonlinear Cases with an Example on Manifold}, we further define $\boxplus$ by the exponential mapping of $\mathrm{SE}(3)$.
Formally these two binary operators are defined as
\begin{equation*}
	\boxminus: \quad  \mathbf{T}_1 \boxminus \mathbf{T}_2 = \mathbf{Log} ( \mathbf{T}_1^{-1} \mathbf{T}2 ), \quad \mathbf{T}_1, \mathbf{T}_2 \in \mathrm{SE}(3)
\end{equation*}
\begin{equation*}
	\boxplus: \quad  \mathbf{T} \boxplus \boldsymbol{\xi} = \mathbf{T} \cdot \mathbf{Exp} ( \boldsymbol{\xi} ), \quad \mathbf{T} \in \mathrm{SE}(3), \boldsymbol{\xi} \in \mathbb{R}^6
	.
\end{equation*}
The constraint ${}^A\mathbf{T}_l = {}^B\mathbf{T}_r$ can be written in the standard form
\begin{equation*}
{}^A\mathbf{T}_l = {}^B\mathbf{T}_r
\Longleftrightarrow
\mathbf{Log}( {}^A\mathbf{T}_l  {}^B\mathbf{T}_r^{-1} ) = \mathbf{0}
.
\end{equation*}
Let us denote $\mathbf{C}_2 ( \mathbf{x} ) = \mathbf{Log}( {}^A\mathbf{T}_l {}^B\mathbf{T}_r^{-1}  )$.
Now we are at the stage to apply the result in Section \ref{section: Extension to Nonlinear Cases with an Example on Manifold} directly.

\textsf{Previous solution and its covariance:}
Without the constraint ${}^A\mathbf{T}_l = {}^B\mathbf{T}_r$, the optimal trajectory for $\mathcal{T}_A$ and $\mathcal{T}_B$ are simply their original forms
\begin{equation*}
\mathcal{T}_A:\quad
   {}^A\mathbf{T}_1^{\star} =   {}^A\tilde{\mathbf{T}}_1
   ,\quad
   {}^A\mathbf{T}_{i,i+1}^{\star} =  {}^A\tilde{\mathbf{T}}_{i,i+1}
   \quad ({i = 1, 2,\cdots N_A})
\end{equation*}
\begin{equation*}
\mathcal{T}_B:\quad
   {}^B\mathbf{T}_1^{\star} =   {}^B\tilde{\mathbf{T}}_1
   ,\quad
   {}^B\mathbf{T}_{j,j+1}^{\star} =  {}^B\tilde{\mathbf{T}}_{j,j+1}
   \quad ({j = 1, 2,\cdots N_B})
   .
\end{equation*}
Let us collect the above result in the variable $\mathbf{x}_1^{\star}$.
To calculate its covariance
$
\mathbb{C}\mathrm{ov}(\mathbf{x}_1^{\star})
$,
we linearize the cost function at $\mathbf{x}_1^{\star}$.
By some trivial linear algebra on Lie group, we get
\begin{equation*}
\begin{aligned}
&\mathbb{C}\mathrm{ov}(\mathbf{x}_1^{\star}) = \mathbf{H}^{-1} \boldsymbol{\Sigma} \mathbf{H}^{-T}
\quad ( \mathrm{by\ (\ref{eq: covariance of least distance optimization})} )
\end{aligned}
\end{equation*}
with
$
\mathbf{H} = \mathbb{I}
$
,
and $\boldsymbol{\Sigma}$ constructed diagonally with corresponding ${}^A\Sigma_i$, ${}^A\Sigma_{i,i+1}$,
${}^B\Sigma_i$, ${}^B\Sigma_{j,j+1}$.
The optimal value $f^{\star} = 0$.

\textsf{New Constraint:}
For the constraint,
let us index the $i$-th component of a trajectory, for example the trajectory $\mathcal{T}_A$ as $\mathcal{T}_A[i]$,
which basically means $\mathcal{T}_A[1] = {}^A\mathbf{T}_1$, $\mathcal{T}_A[2] = {}^A\mathbf{T}_{1,2}$, $\mathcal{T}_A[3] = {}^A\mathbf{T}_{2,3}$ and so forth.
Let us further denote
\begin{equation*}
\mathbf{C}_2 (\mathbf{x}_1^{\star}) = \mathbf{Log} ( {}^A\tilde{\mathbf{T}}_l {}^B\tilde{\mathbf{T}}_r^{-1} ) = \eta
.
\end{equation*}
The linearization of $\mathbf{C}_2 (\mathbf{x}_1)$ at $\mathbf{x}_1^{\star}$ consists of trivial linear algebras based on the BCH formula and the adjoint operation of $\mathrm{SE}(3)$ \cite{chirikjian2011stochastic}\cite{barfoot2017state}.
We provide the analytical Jacobian here and ignore the calculation details for simplicity:
\begin{equation*}
	\mathbf{A}_2 (\mathcal{T}_A [i]) = - \mathbf{J}_l(\eta)^{-1} \mathbf{Ad}({}^A\mathbf{T}_i^{\star})
	,\  (i = 1, 2,\cdots l)
\end{equation*}
\begin{equation*}
\mathbf{A}_2 (\mathcal{T}_B [j]) = \mathbf{J}_l(\eta)^{-1} \mathbf{Ad}
	({}^A\tilde{\mathbf{T}}_l {}^B\tilde{\mathbf{T}}_r^{-1} {}^B\mathbf{T}_{r - j + 1}^{\star})  )
	,\  (j = 1, 2,\cdots r)
\end{equation*}
with $\mathbf{J}_l(\cdot)$, and $\mathbf{Ad}(\cdot)$ being the left-hand Jacobian and the adjoint operation of $\mathrm{SE}(3)$ respectively.

\textsf{Apply Eq. (\ref{equation: predicting objective function change in nonlinear least distance optimization}):}
Finally, the change of optimal values $\Delta f$ can be calculated by Eq. (\ref{equation: predicting objective function change in nonlinear least distance optimization}),
and the optimal value of the problem (\ref{eq. optimization problem: cost of aliging trajectory gaps}) can be approximated by
$
f^{\star\star} \approx f^{\star} + \Delta f = \Delta f
$.
Note that for this specific example, we do not need to solve any optimization problem. The only requirement is to linearize the constraints,
and then apply (\ref{equation: predicting objective function change in nonlinear least distance optimization})!

\section{Numerical Experiments}
\label{section: Numerical Experiments}

In this section, we provide numerical experiments for the trajectory alignment problem defined in (\ref{eq. optimization problem: cost of aliging trajectory gaps}).
We use Matlab to simulate two random robot trajectories, by firstly rotating the robot by a random angle, and then moving the robot $1m$ forward.
These two trajectories are intersected in the middle in case they are completely apart from one other.
Afterwards, we add Gaussian additive noise with standard deviation $0.1m$ on the translational part, and $0.01rads$ on  the rotational part.
We use an Intel i5-5300U CPU @ 2.30GHz $\times$ 4,
with Ubuntu 16.04 LTS distribution
to run all the experiments.

\subsection{An Intuitive Example: Trajectory Change, Predicted Cost and Real Cost}

First, let us consider a simple example where each trajectory contains $20$ poses.
We choose several pose combinations to demonstrate the idea as shown in Fig. \ref{lb: Examples of aligning two trajectories, trajectory change and cost}.
We visualize the change of the trajectory with respect to the cost predicted by the metric,
as well as the real cost obtained by completely solving the problem (\ref{eq. optimization problem: cost of aliging trajectory gaps}).
From Fig. \ref{lb: Examples of aligning two trajectories, trajectory change and cost},
we would say that the metric is quite robust in general;
even for the cases of misalignments, like pose $1$ and $20$, the resulting trajectories are completely twisted, while the metric is still pretty accurate.
Note that the initial noise free trajectories are perfectly aligned at the $10$-th poses,
where the problem attains the lowest cost (with noise) and the metric is particularly precise.

\begin{figure}[t]
	\includegraphics[width=.45\textwidth]{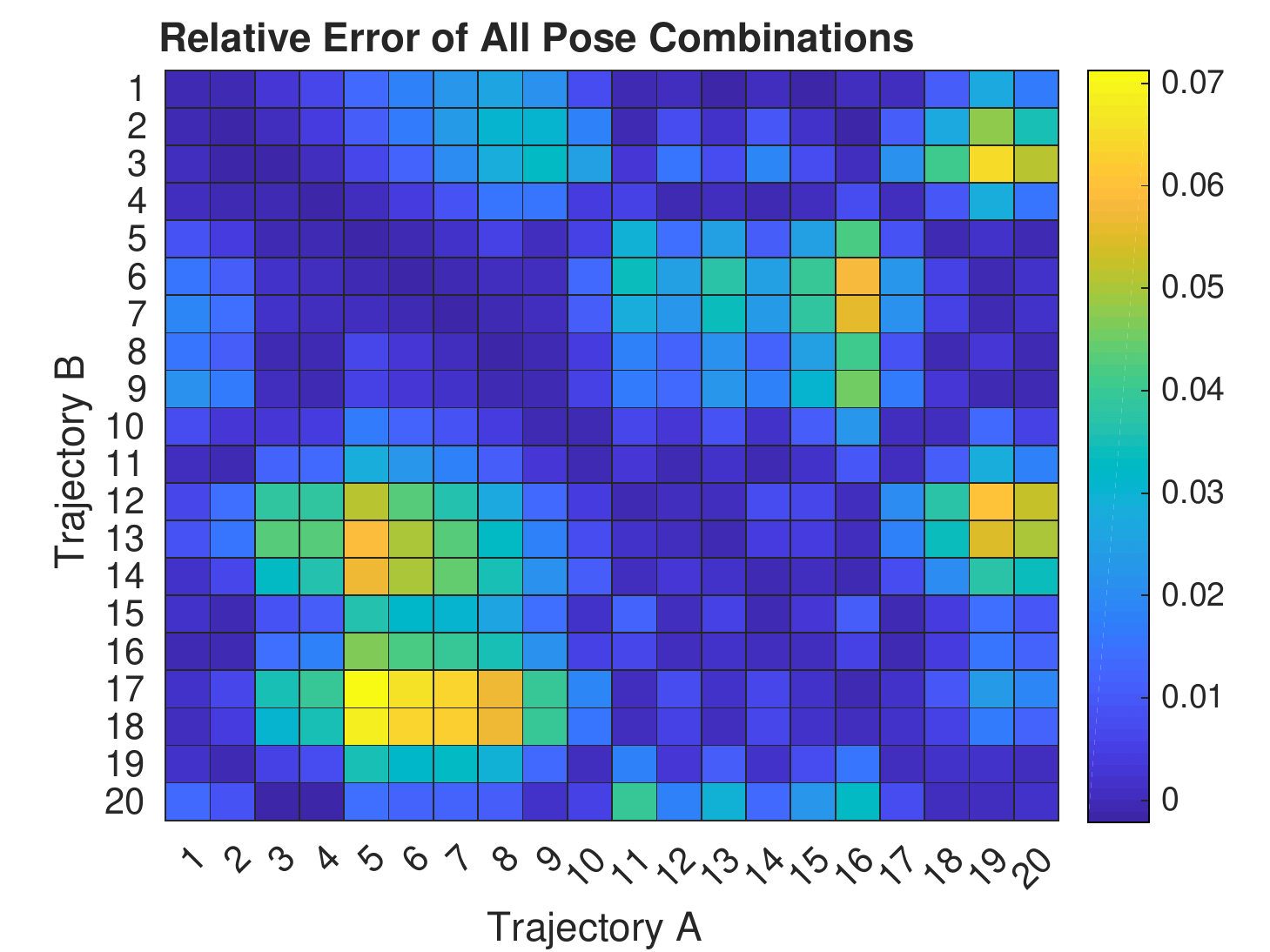}
	\caption{The relative error of the predicted cost in each pose combination.}
	\label{The relative error of the predicted cost in each pose combination.}
\end{figure}

Then we calculate the cost for all possible pose combinations, which is $20 \times 20$, and report the corresponding relative error in Fig \ref{The relative error of the predicted cost in each pose combination.}.
It can be envisioned that there exists many unrealistic alignments within these $20 \times 20$ combinations, however the relative error is essentially very low, with a maximum around $10 \%$, which is sufficient for any further decision makings.
At last, with the Matlab code, it takes only $3$ seconds to compute the metric for all the $20 \times 20$ cases, while the corresponding time consumed by solving the problem is $26$ seconds.

\subsection{Computational Time and Relative Error with Different Trajectory Lengths}

To quantify the effect of trajectory lengths on the accuracy of the metric,
we simulate several trajectories with $20$, $50$, $100$ and $200$ poses respectively.
As before, we compute the alignment for all possible pose combinations for each case.
The relative error between the predicted cost (by the metric) and the real cost (by solving the problem) is reported in Fig. \ref{lb: trajectory alignment. relative error of the metric, with respect to different lengths of trajectories.}.
It can be seen from Fig. \ref{lb: trajectory alignment. relative error of the metric, with respect to different lengths of trajectories.} that for most of the cases, the relative error is still reasonable. However as the length of the trajectory grows, the ``unrealistic alignments" (i.e. misalignments) can lead to a substantially amount of outliers in the boxplot.
An explanation to this phenomenon is that because the metric is essentially built on the result of linear cases, if the alignment is too far, the linearization point changes dramatically and the linear result is less accurate.
However, in general, the predicted cost still provide a reasonable approximation to the real cost.

The computational time is reported in Table \ref{lb: trajectory alignment, the computational time of each trajectory lengths.}.
For each case, we show the overall computational time for all possible pose combinations, as well as the averaged time for obtaining a cost for a single pose combination.
Overall, calculating the cost by the metric is at least one magnitude faster than by solving the problem directly.
This timing saved here can be extremely beneficially for large problem instances.

\begin{figure}[t]
	\includegraphics[width=.415\textwidth]{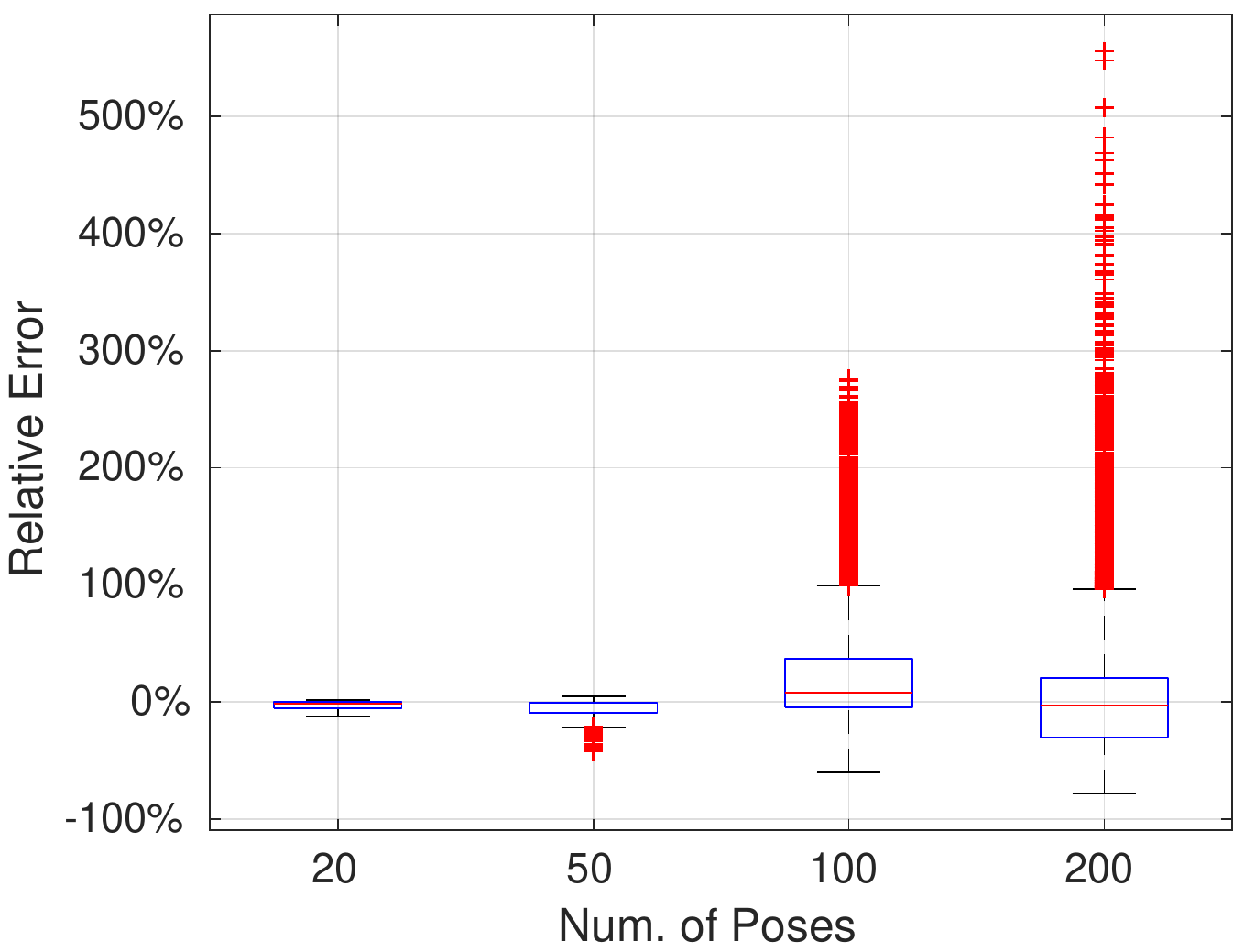}
	\caption{The relative error of the metric with respect to different lengths of trajectories. Note here we compute all possible pose combinations, for example the number of possible alignments for $200$ poses are $40000$. This may of course include some ridiculous combinations, which attributes to the outliers (marked in red) in the plot.}
	\label{lb: trajectory alignment. relative error of the metric, with respect to different lengths of trajectories.}
\end{figure}

\renewcommand{\arraystretch}{1.15}
\begin{table}[t]
	\caption{The Computational Time (Unit: Seconds)}
	\label{lb: trajectory alignment, the computational time of each trajectory lengths.}
	\centering
	\begin{tabular}{ccccc}
	\hline\hline
	Num. of Poses & 20   &   50   &   100   &   200  \\
	\hline
	By metric overall & 2.49116   &   36.7919   &   294.033   &   2457.14  \\
	By solving overall  &  22.5889   &   460.228   &   4799.49   &   42428.2 \\
	\hline
	By metric avg. & 0.00623  &  0.0147  &  0.0294  &  0.0614 \\
	By solving avg. & 0.0565  &   0.184   &   0.480 &   1.06 \\
	\hline
	\hline
	\end{tabular}
\end{table}
\renewcommand{\arraystretch}{1.0}

\section{Discussion and Conclusion}

In incremental scenarios, it is also possible to compute an approximate solution to the second subproblem
(i.e., Eq. (\ref{minimal norm optimzation: full stacked measurement})
(\ref{Eq: least distance optimization, full measurement matrix})),
exploiting special structure and sparse linear algebraic techniques
\cite{gill1974methods}\cite{cassioli2013incremental}\cite{kaess2008isam}\cite{kaess2012isam2}\cite{polok2013incremental}.
While this idea works and is a general practice,
the main insight of this paper over previous publications is that: It is actually possible to derive an analytical equation (in closed form) which quantifies the change of optimal values directly.
Besides, the proposed method is more advantageous in face of many new measurement candidates,
which can share
$
\mathbf{x}_{1}^{\star}
$
and
$
\mathbb{C}\mathrm{ov} (\mathbf{x}_{1}^{\star} )
$.

In the proposed equations,
the computational bottleneck of $\Delta f$ lies in the calculation of
$
\mathbb{C}\mathrm{ov} (\mathbf{x}_{1}^{\star} )
$,
which can be mitigated by computing a marginal covariance \cite{ila2015fast}\cite{bai2018predicting},
i.e., only computing
$
\mathbf{A}_2
\mathbb{C}\mathrm{ov} (\mathbf{x}_{1}^{\star} )
\mathbf{A}_2^T
$
explicitly.
It is also of great interests to further compare the numerical procedure to calculate an updated solution to the second subproblem,
with that to calculate
$
\mathbf{A}_2
\mathbb{C}\mathrm{ov} (\mathbf{x}_{1}^{\star} )
\mathbf{A}_2^T
$.
Further research in this direction will help understanding the connection between updating the solution and 
updating the optimal value.

In nonlinear cases, the equation is derived by linear approximation.
The numerical experiment shows good accuracy in general,
but still there are many failure cases as shown in Fig. \ref{lb: trajectory alignment. relative error of the metric, with respect to different lengths of trajectories.}.
One explanation is to use the change of linearization points.
However, Fig. \ref{lb: Examples of aligning two trajectories, trajectory change and cost} shows that the approximated $\Delta f$ is actually quite tolerant to linearization changes (see the twisted trajectories).
Therefore
it is reasonable to assume that there exists a more sophisticated mechanism that has not been well-understood yet.
A theoretical investigation on the accuracy bound of $\Delta f$ is still required.

To conclude, this paper provides the closed form of the change of optimal values in the incremental minimum norm optimization and least distance optimization in the linear case.
The change of optimal values in nonlinear cases are approximated via linearization. 
These results yield the possibility of using the change of optimal values as a pre-calculated metric, in incremental (or online) applications.


\end{document}